**Christian Retoré**
Université de Bordeaux —LaBRI CNRS
351, cours de la Libération 33405 Talence cedex
**christian.retore@labri.fr**

**Variable types for meaning assembly:
a logical syntax for generic noun phrases introduced by *most***

### 1. Overview

This paper proposes a way to compute the meanings associated with sentences with generic noun phrases corresponding to the generalized quantifier *most*. We call these generics *specimens* and they resemble stereotypes or prototypes in lexical semantics. The *meanings* are viewed as logical formulae that can thereafter be interpreted in your favourite models.

To do so, we depart significantly from the dominant Fregean view with a single untyped universe. Indeed, our proposal adopts type theory with some hints from Hilbert ε-calculus (Hilbert, 1922; Avigad and Zach, 2008) and from medieval philosophy, see e.g. de Libera (1993, 1996). Our type theoretic analysis bears some resemblance with ongoing work in lexical semantics (Asher 2011; Bassac et al. 2010; Moot, Prévot and Retoré 2011).

Our model also applies to classical examples involving a class, or a generic element of this class, which is not uttered but provided by the context. An outcome of this study is that, in the minimalism-contextualism debate, see Conrad (2011), if one adopts a type theoretical view, terms encode the purely semantic meaning component while their typing is pragmatically determined.

## 2. Type theory as a syntax of semantics.

In the viewpoint of Montague (1970), a syntactic parse tree is turned into a reference, that is a set of truth conditions, or an interpretation in possible worlds, via an intermediate step, a logical formula written as a lambda term, that one could call a *logical form.* In Montague's analysis, intermediate steps ought to be wiped off once the computational process yielding references is over. Our article is on the contrary focused on this logical form and how it can be computed — this step is reminded in Moot and Retoré (2012). What should be the logical formulae depicting the meaning? There are several reasons for addressing this question:

- Firstly, given the subtlety of the questions we wish to consider, involving generalised quantification and vagueness, the right notion of reference, if any, is going to be tricky.
- For such issues, it likely that other interpretations, not in terms of truth and reference, but in terms of argumentation are better designed — for instance, an interactive view of semantics has been proposed in Lecomte and Quatrini (2011) using *ludics*: roughly speaking, formulae are interpreted as sets of possible proofs and refutations.
- In order to properly depict an interpretation, be it a referential or an interactive one, we first need to know precisely which formulae we wish to interpret.

The logical syntax of semantics is usually split into two systems:

1. The logic in which we express the logical form, the meaning of the sentence, discourse, etc. In Montague semantics it is usually higher order predicate logic — although some use first order logic via a reification of predicates. This logic is used as a language, which can be interpreted in models.
2. The glue logic, also known as metalogic, is usually presented as a lambda calculus — simply typed lambda calculus, in the case of Montague. Via the Curry-Howard isomorphism, well presented in Girard *et al.* (1988), these terms are formal proofs in propositional intuitionistic logic with two base types: **e** for entities and **t** for propositions. This logic is used as a calculus for meaning assembly, since such a view is mainly compositional: the meaning of the whole is a function of the meaning of its parts and of the syntax — although context should also be taken into account, as we shall see.

We adopt a similar viewpoint, but there will be a difference: while there is no interplay between the two logical layers in standard Montague semantics, we shall have some interaction, because the metalogic is richer:

1. The logical forms take place into multisorted higher order logic.
2. The metalogic, the lambda calculus for meaning assembly will be second order propositional intuitionistic logic, that is second order lambda calculus with many primitive types for entities (we give precise definitions of this system in the next section).

The need for having many primitive types is guided by the following critique of Frege's single sorted universe. We think that language rarely or even never quantifies over all entities, but over entities of a given type. Frege used a trick to reformulate without types what natural language or even mathematics expresses as $\forall a{:}A\ P(a)$ with A a type. He writes this as $\forall x\ [A(x) \Rightarrow P(x)]$. But this formulation is *ad hoc*, and it does not extend to generalised quantification, since the two following propositions, in which $\angle$ means "most" are not equivalent:

- $\angle a{:}A\ P(a)$         *Most students go out.*
- $\angle x\ [A(x) \Rightarrow P(x)]$ *For most entities, if they are students they go out.*

This is a reason why the standard approach view generalised quantification as a function of two predicates, the restrictor and the predicate itself, see e.g. Peters and Westerstahl (2008), Szabolcsi (2010), Mari (2011). Thus the study of generalised quantifiers cannot treat quantification as the usual logical setting does, with quantification acting on a single predicate.

Our preference for a rather natural quantification over a type (a class, a set) rather then for an absolutely universal one is reminiscent of medieval philosophy. Indeed, philosophers around Avicenna, in particular his student Abu'l-Barakat al-Baghdadi, said that properties should always be asserted of an object as being a member of some class (perhaps nowadays they would say *type*), and they draw a distinction between homogenous predicates that apply to a class and its subclasses and predicates that apply across classes (De Libera, 1993). In comparison with the above formulation by Frege, there is quite a difference: not any formula with a single variable defines a type, types are much more restricted, and so are comparison classes.

If we adopt a type system where not any set is a type, then we shall probably want some flexibility, enabling an object of a given type to be considered as a member of another type when the context requires such a type change. Such transformations require operations on families of types, hence a second order type system. Our idea behind this framework is to draw a border between *semantics* and *pragmatics:* lambda terms compute *semantic* representations and types are irrelevant for computing beta reduction and substitutions. The types are used to filter impossible compositions (like "*The table shouts*"), and they are often determined the context.

Our notion of lexicon introduced in Bassac et al. (2010), and extended by Moot and Retoré (2011b) is a rather simple extension of a Montagovian lexicon. In addition to the standard lambda term expressing the argumental structure of a word, each word is endowed with a finite set of lambda terms that allow type transformation in case of solvable type mismatch. As an example, let us simply say we provide a correct account of sentences like:

(1) Liverpool is a poor town and an important harbour.
(2) * Liverpool defeated Chelsea and is an important harbour.

Base types include *town*, *location*, *people*, *soccer_team*,… and many more, since any class that rather naturally comes to mind should be a base type. The type *town* can be converted to any other type in this list, but some transformations are declared as exclusive by the lexicon. For instance the transformation of a town into a *soccer_team* is incompatible with the other transformation that a town may undergo. [We do not handle examples involving phenomena like contrast "*The small town of Manning defeated the much larger urban school of Davenport by a score of 43-46.*" (http://en.wikipedia.org/wiki/Manning,_Iowa) nor strong contextualisation "*Barcelona won four UEFA champions league and organised the Summer Olympic Games of 1992*".]

Our model was initially designed to handle meaning transfers in compositional lexical semantics and thereafter applied to other questions like plurals, or here, quantification and generic NPs. By generic NPs we mean expressions that correspond to an implicit or explicit "most of" quantifier:

(3) *Most* Lords are appointed by the Queen. (www.parliament.uk)
(4) *The* Brits love France. (www.**brits**-in-**france**.net) [Here we shall only consider the reading "*Most Brits love France.*" while this sentence admits other readings: "*Brits prefer France to other countries*" or "*Brits like France more than citizens from other countries do.*", and possibly others.]

Our idea is to consider a generic element corresponding to *most*, as Hilbert $\varepsilon x. P(x)$ and $\tau x. P(x)$ have used to model existential $\exists$ and universal quantification $\forall$ (Hilbert 1922, Avigad and Zach, 2008). The universal generic $\tau x. P(x)$ is an element such that $P(\tau x. P(x))$ holds whenever every x has the property P. Symmetrically, $P(\varepsilon x. P(x))$ holds whenever there exists an *x* enjoying the property P. In linguistics, the closest construct is the $\iota$ of von Heusinger (2007), which combines Hilbert's ε with a choice function: it is used for picking an element among those that have a given property, when there are some. Demonstratives and definite articles behave like this: "*Give me the red pen.*" Context and pragmatic principles have to be used to know which element is selected: this choice is left out of the purely semantic calculus.

## 3. Metalogic: the system F of variable types

Types are defined from the base type **t** and many base types for entities, since semantics is filtered by type mismatch. A typical mismatch occurs when syntax determines application of a predicate ranging over type $e_i$ objects to type $e_j$ objects. We thereafter need some uniform operations applying to all types, and the type of such an operation is *"for all type X T[X]"* : such a universal operation can be specialized to any type X. For instance, the aforementioned ι of von Heusinger could receive the second order type *"for all type X. X"* that is, for all type X the operator ι yields an object of type X — see the end of the section for a discussion. The type *"for all type X. T[X]"* is written *"ΠX. T[X]"* — we do not write *"∀X. T[X]"* because we also need the usual symbol "∀" for quantification in the logical form. This second-order lambda-calculus called *system F,* was introduced by Girard (1971), but a good and freely available reference is Girard *et al.* (1988).

Types are defined as the elements of the smallest set containing:
- Base types, that are **t** and $e_j$ (a finite but large set of entity classes, which can be viewed as some linguistic ontology)
- Type variables (variables ranging over types) denoted by upper case letters $Xi$ (a countable number of them is required, as usual with variables)

The set of types is closed under the following operations:
- Whenever *T1* and *T2* are types, *T1→T2* is a type as well. This type represents functions that map *T1* objects into *T2* objects. The type that we write *T1→T2* is sometimes written *<T1,T2>* in the Montagovian literature.
- Whenever *X* is a type variable and *T* a type, *ΠX. T* is a type — usually *T* depends on X.

The usual Montagovian system uses simply-typed lambda-calculus, which, via the Curry-Howard isomorphism — see e.g. Girard *et al.* (1988) — corresponds to intuitionistic propositional logic with two base propositions, **e** and **t** and implication → as the only connective. In our model, one may quantify over propositional variables. The corresponding logic is intuitionistic second order propositional logic with implication and universal quantification in which other connectives can be encoded such as existential quantification, conjunction and disjunction.

The terms typed within this type system can be viewed, by the Curry-Howard isomorphism, as proofs in intuitionistic second order propositional logic. They are defined as follows where both expressions *"t:U"* and *"$t^U$"* denote a term *"t"* of type *"U"*.

Firstly, terms include variables and constants:
- For every type, we have a finite and possibly empty set of constants:

- - Some of them correspond to constants of the semantic language for logical forms (proper names, predicates, higher order predicates, …). For instance, one could have a constant *love* of type *Animated* → ***e*** → ***t***
  - Some of them correspond to logical operations (connectives, quantifiers etc.).
- For every type we have a countable set of variables of this type.

The set of terms is the smallest set containing the variables and constants above that is closed under the following operations:

- Whenever *f* is a term of type *T→U* and *t* a term of type *T*, the expression *(f t)* which reads "*f applied to t*" is a term of type *T*.
- Whenever *u* is a term of type *U* and *x* a variable of type *T*, the expression *λx. u* is a term of type *T→U*. It denotes "*the function that maps x to u*".
- Whenever *t* is a term of type *ΠX. T* and *A* is a type, the expression *t{A}*, "*t applied to the type A*", or "*t specialised to the type A*", is a term of type *T[X:=A]*, that is the type *T* in which any occurrence of the type variable *X* is replaced with the type *A*.
- Whenever *t* is a term of type *T* and *X* a type variable not appearing in the type of any free variable in *t*, the expression *ΛX. t* is a term of type *ΠX. T* "*the generalisation of t (in X)*".

This definition warrants a few comments.

- A term *u* with a universally quantified type *ΠX. T* may be viewed as a generic term, or a universal term, which can be specialised to a type *A*, by application of *u* to a type *A*, written *u{A}*.
- As one can see, there is a restriction in the generalisation rule. This restriction means that, if in *t* we are assuming some properties of *X*, then we cannot generalise *t* to all types *X*: the presence of a free variable whose type involves *X* precisely expresses that we are assuming that some specific properties of the type X hold.
- Observe that the rules for specialisation and generalisation look a lot like the usual lambda calculus rules for application and abstraction.

If one thinks of types as propositions and of terms of type *P* as proofs of *P*, as the Curry-Howard isomorphism allows, to which deductive system does this lambda calculus correspond? We have propositional constants, and propositional variables and the only connective is the implication "→" together with the universal quantifier over propositional variables "Π": specialisation and generalisation are respectively the quantifier elimination and introduction rules. From these two connectives, one can define conjunction, disjunction, and existential quantification over propositional variables, see e.g. Girard *et al.* (1988). Since there is no classical axiom like *tertium non datur,* it is an intuitionistic sytem. Hence, while the metalogic behind standard Montague semantics is propositional intuitionistic logic with arrow as the only connective, here the meta logic is quantified propositional intuitionistic logic — in which the other connectives and second order existential quantification can be encoded.

As in the usual Montagovian setting, only closed normal terms of type *t* correspond to semantic representations, which are logical formulae. Here as well, we need some reduction rules which compute semantic representations, i.e. formulae, by substituting types for type-variables, and terms for variables. The two reduction rules are very similar: one of them is standard beta reduction and the other is more or less the same, but on *type variables* and *types*.

1. $((\lambda x^A. t^B) u^A)$ reduces to $t[x:=u]$, that is t in which any occurrence of *x* is replaced by the term *u*. Observe that necessarily *x* and *u* are of the same type *A*, and that the type *B* of the whole term is preserved under reduction.
2. $(\Lambda X.t^T)\{A\}$ reduces to $t[X:=A]$, that is *t* in which any occurrence of the type variable *X* is replaced with the type *A*. Observe that the type of the term, $T[X:=A]$, is preserved under reduction.

System F may seem a bit unsafe at first glance. Indeed we are defining the type $\Pi X.T[X]$ from all the types $T[X]$ with X ranging over all types, including $X=\Pi X.T[X]$ which yields $T[\Pi X.T[X]]$, and one may fear this impredicativity. However system F does not collapse. A first argument is that it is possible to construct a model, coherence spaces - see e.g. *Girard et al.* (1988) - with types as structured sets and proofs/terms as particular objects of these structured sets. Another argument shows that there is no collapse and will be useful to us thereafter: the terms of system F enjoy a strong normalisation, a confluent one yielding a unique normal form however one proceeds (Girard, 1971). This shows system F is safe: if there were a collapse, or proof of something false, then there would be a normal one, and given the shape of normal terms/proofs it is easily seen that there cannot exists such a term/proof.

## 4. The lexicon and the logic of semantic representations

The metalogic is second order propositional intuitionistic logic a.k.a. second order lambda calculus, but in which logic are we going to express the meaning of sentences? Actually, we do not depart much from the standard Montagovian approach. We use higher order logic but a multisorted one: indeed, we have several base types for entities with relations among them. These relations between types are encoded by functions mapping one type to another and they represent meaning transfers that are reminiscent of Nunberg (1995). Here as well, both the argument and the predicate can contribute to the transfer of meaning, but the main difference is that we integrate these meaning transfers into a broader compositional type-theoretic framework, fully formalized and with an explicit computational mechanism — this lead us to some schematic model, leaving out some semantic subtleties.

An entry in the lexicon associates to each word a main lambda-term, that of standard Montague semantics, within a richer type system. But there might be other lambda terms for turning objects of a given type into the same object considered as a member of some other type. For instance a book is of type *Book*, but it can be turned into a material thing that can be heavy, or into contents that may be interesting. Ontological inclusions like *Cars* into *Vehicles* are also encoded that way. In order to block infelicitous compositions like "*The table barks*", we have many base types, and the argument of "*barks*" should be of type *"Dogs"*. But these type constraints should also be relaxed to allow composition of sentences like "*I am parked in front of building 20*" — both the predicate/function and the argument can provide the transfer. Finally, some meaning transfers are irreversible, and block other possible transformations: once a *town* has been understood as a *soccer team*, it cannot be considered as a *location* or as *its major*.

In order to illustrate the kind of operations such a system is able to model, let us see how we can conjoin two predicates that apply to different kinds of objects, like *heavy* and *interesting* which respectively apply to *material things* (M) and to *abstract contents* (A). This conjunction can only be considered if an object can be viewed both as a material thing and as an abstract content, as *books* (B) can. Given

- two predicates, *h* of type M→**t** and *i* of type A→**t**,
- two maps, *m* from a type B to M and *a* from B to A,
- *b* of type B

the conjunct can be expressed as *(h (m b)) and (i (a b))* or rather, using lambda calculus prefix notation, as

$$(and\ (h\ (m\ b)))\ (i\ (a\ b))$$

In this expression, "*and*" of type $t \to (t \to t)$ is the standard conjunction of two propositions. As we are able to do so in any such situation, the lambda term AND for the word "and" in a sentence like

(5) This book is heavy and interesting.

is

$$\text{AND: } \Pi A\ \Pi M\ \lambda i^{A \to t} \lambda h^{M \to t}\ \Pi B\ \lambda b^B\ \lambda a^{B \to A} \lambda m^{B \to M}$$
$$(and^{t \to (t \to t)}\ (h\ (m\ b)))\ (i\ (a\ b))$$

Observe that the variables whose type contains a type variable *X* are bound before quantifying over the type-variable X. The functions *a* and *m* from the type B (*book* in our example) to, respectively, the types *A* and *M* are provided by the entry *book* in the lexicon. With the strength of second order quantification, a single term is enough, since this single term can be specialised to the types of the situation under consideration. For our book example, AND should be successively applied to the following five terms and three types in this order:

AND {Abstract}{Material}

$$\text{interesting}^{\text{Abstract}\to t} \text{heavy}^{\text{Material}\to t}$$
$$\{\text{Book}\}$$
$$\text{this\_book}^{\text{Book}} \text{to\_contents}^{\text{Book}\to A} \text{to\_material}^{\text{Book}\to M}$$

Another example, for the syntax of semantics, is the modelling of the universal quantifier ∀. In an unsorted logic, one only needs a single quantifier for individuals, which from the viewpoint of lambda calculus is a constant of type *((e→t)→t)* and others for higher order quantifications, e.g. *(((e→(e→t))→t)→t)* for quantifying over transitive verbs. Using system F, a single constant for universal quantification is needed. A single constant ∀ of type *ΠX. (X→t)→t* is enough since it can be specialised to a type *U* if one wants to quantify over *U*-objects by specialisation to *U*. Indeed, the term *(∀ {U})* is of type *(U→t)→t* and *(∀ {U})* represents universal quantification over U-objects. The same consideration applies to von Heusinger's ι. We said its type should be *ΠX. X* but we should also have a constant ι receive the type *ΠX. (X→t)→X*, which when applied to a property of *X* yields an element of type *X* enjoying the property *P*.

## 5. Some questions on this computational model of semantics

This reorganisation of compositional semantics and of lexical semantics deserves some comments. Firstly one should really make a distinction between the metalogic, quantified propositional intuitionistic logic, where the proofs/terms are relevant objects, and the logic of semantic representations, higher order multisorted logic, where only the formulae are relevant. But, in our model, there is an unusual interplay between them: for any judgement "u *is of type P*" i.e. *u:P* one can consider a predicate *P* and state as an axiom the formula *P(u)*. On the other hand, not any formula with a single variable defines a type: types are much more constrained, and from an intuitive viewpoint, they should correspond to cognitively accessible classes.

A frequently asked question about our model is whether it has something to say about subtyping. Indeed, it is known that, despite some attempts by Cardelli *et al.* (1994) and by Soloviev and Luo (2000), "subtyping" does not get along well with second order typing. Firstly, one should draw a strong difference between the semantic notion of subtyping and the technical notion of "subtyping" in typed functional programming — remember that system F and similar lambda calculi can be viewed as functional programming languages, so the confusion is possible. Subtyping in functional programming is supposed to fit in with the functional types: the subtypes of *A→B* should be inferred form the subtypes of *A* and those of *B*. Do the linguistic subtypes of eating verbs *food → human → t* derive from the subtypes of *food* and of *human beings*? We think that manners are more relevant for classifying eating verbs than the nature of their subjects and objects. Furthermore, it is likely that the linguistic IS_A relation is idiosyncratic and much more constrained than the real world ontological relations.

There exist alternative solutions to our proposal for handling both compositional and lexical semantics and pragmatics: the one by Asher (2011), using categorical logic, and the one by Luo (2011) using type theory. Although such type theories are weaker as logical systems, because they do not allow quantification over any type, these type theories have many structure-building rules, many reduction rules, and many variants. The structures they offer for linguistic modelling can be encoded in system F, except *dependent types*: if the only dependent types one needs are just records, they are already present in system F, but if real dependent types are absolutely needed — this is questionable — dependent types can safely be added to system F, since all of this is included in the calculus of constructions of Coquand and Huet (1988)— see e.g. Bertot and Castéran (2004).

The system F also raises complexity issues, in particular for typing a given pure lambda term, or because it contains functions requiring an exponential number of reductions to compute their results. Our model does not have to face these issues. Indeed, we only reduce terms that are obtained by inserting lambda terms from the lexicon into a syntactic tree: there are neither problematic functions in the lexicon, nor in the syntactic tree/term, and we never try to compute the type of an untyped lambda term. There can be several choices of meaning transfers in case of type mismatch, but this is just the unavoidable syntactic, semantic and pragmatic ambiguity of human languages.

## 6. A logical syntax for generics introduced by "most"

Although there is not a clear cut-off between the two constructions, following English usage see e.g. the Grammar quizzes by Sevastopoulos (2012), we make a distinction between noun phrases introduced by "most" and those introduced by "most of". Indeed, while the class associated with the "most" object is natural and immediately apprehended, the class associated with the "most of" object may be a complex one, defined by a formula.

(7) Most students passed logic.
(8) Most of the students that passed logic passed algebra.
(9) Most of the students went to the university party.

Observing (7) and (8) one has the feeling that the "most of the" construction may apply to any class defined by a complex property of one entity, while the plain "most" applies to a natural class. This impression is confirmed by comparing (7) and (9): in (9) one has the impression that word "student" is used as a property of entities of a wider class, e.g. a class including students, professors, administration… The distinction will be reflected in our system by the distinction between types and formulae with a single free variable.

### 6.1. The logical syntax of bare "most" generics

Our idea is that "most" always refers to a type and not to any set that could be defined by a formula with a free variable using the comprehension scheme: this is the "most_of_the" quantification. Consequently, combining our treatment of quantification, a single quantifier which can specialised to any type, and generic elements like the τ and ε of Hilbert (1922) and the ι of von Heusinger (2007), we propose to use a constant "∠", read as *specimen of*, of type $\Pi X.\ X$. As von Heusinger's ι whenever it is applied to a type $A$ it yields the *specimen* in $A$: $∠\{A\}$ is of type $A$. Using predicate logic rather than type theory, if A were a property, that is a unary predicate, this would be written ∠x. A(x) in the style of von Heusinger's ι.

As opposed to the work on generalized quantifiers (see Keenan, Peters and Westerstahl) the generalised quantifier is defined from a single predicate/type A, and not as a function of two predicates: we only use the first of them, which is assumed to be a type in the bare "most" case.

Let us come back to the standard examples and compute the readings. A syntactic analysis, e.g. a categorial one, of example (4) "*The Brits love France*." will yield the linear lambda term: *((loves France) (the Brits))* — syntactically the verb is first applied to its object and then to its subject. Firstly, observe that despite the quantifier, there is no need for type raising. Indeed, *the Brits* will be a virtual element of type *Brit*, which is a subset of *Human_beings* which is itself a subset of *Animals* which are the right class for subjects of the binary predicate *love*. There is no need to restrain the object of the binary predicate *love*. Hence the lexicon, within the entry *Brits*, provides two morphisms *h* of type *Brits → Humans* and *a* of type *Human → Animals*. These transformations, which are just type inclusions encoding the subtyping relations, are quite particular and transparent: hence it is possible to have a general rule saying that whenever a function applies to a type then it applies to any type included in it, without writing the transformation(s), and this amounts to having as many functions "*love*" as there are subtypes of its arguments. Thus the term denoting the semantic representation before reduction is either

*(λy λx loves x y) France (a(h(∠{Brits})))*

which is well typed, and which reduce to the lambda term

*(loves (a(h(∠{Brits} )))) France*

— with implicit type inclusion the first term would be *(λy λx loves x y) (∠{Brits}) France*, and the reduced one would be *(loves (∠{Brits} )) France*.

The constant "∠" is the lambda term associated with *the,* meaning *most of* and this term produces the specimen associated with *the Brits* — it does not prevent *the* from having other behaviours meaning a definite set,

"*all the*"… etc. Other examples introduced by "*most*" are handled exactly the same way. Now let us turn our attention to generics introduced by "*most of the*" with an NP thereafter.

## 6.2. The logical syntax of "most of the" generics

For processing an example like (8), the model is quite similar, but the constant ∠ to be used has a slightly different type: ∠: ΠX. (X→t)→X. This constant ∠ takes a property of *X*-objects, and returns the specimen of the corresponding subset of *X*. Assume a categorial analysis with words replaced by their semantic lambda-terms yields:

$$passed(\angle(\lambda x{:}student\ passed(x,logic)),topology).$$

Letting *s* be a shorthand for the "*most of the*" generic i.e. s=∠*(λx:student passed(x,logic))*, the typing ensures that s is of type *student*, but it is a bit tricky, in plain lambda calculus, to also get the information *passed(s,logic)* and to produce, as final semantic representation, *passed(s,logic) and passed(s,topology)*. It would be much more convenient in the lambda DRT of Muskens (1996). This is the way semantics is implemented in the categorial parser of Moot, see e.g Moot *et al. (2011)*, but an explanation of lambda DRT would be too lengthy to be included in this paper.

## 7. A word on interpreting the generic element

Although we said earlier that we do not yet have any proper interpretation of the *specimen* of a type, let us say a word about the reference of this generic element, that is its truth conditions, and thereafter on its possible interpretation in interactive terms, that is its proof theoretical usage.

### 7.1. When is "most of the A are B" true?

Firstly, as opposed to some of the literature, but in accordance with Solt (2009), we strongly assert that "*most*" is much more than "*the majority of the*". For instance, after an election won with 53% against 47% of the vote, one cannot say that most of the electors voted for the first winner. "*Most*" is a vague quantifier and starts to be true from a percentage that varies according to the class and the predicate. Possibly "70%" is enough to say that *most students contracted the flu*, but at the same time "70%" is not enough to assert that most of the students passed the exam.

Secondly, as opposed to what is said by most of (!) the literature, it has little to do with cardinality but rather concerns measure, and our apprehension of the class, which can be infinite. For instance, one can find, even in advanced maths books the statement "*most numbers are not prime*" (*An invitation to modern number theory*, by Steven J. Miller and Ramin

Takloo-Bighash). What does that mean? It is known from the Ancient Greek mathematicians that there are as many numbers as there are prime numbers. We need to consider a measure for the whole class and consider that the measure of the relevant subset is a large percentage of the measure of the whole or to consider that the limit of some cognitive perception of the class: the statement on prime numbers simply means that 0 is the limit of the proportion of prime numbers between 0 and $n$ when $n$ approaches infinity. A natural notion of measure is the usual mathematical notion of measure used in probability theory. It also applies to infinite sets, and makes it possible to consider limits, to have subsets with the same cardinal but with different measures,... as natural language does.

Given our choice to handle "*most of*" quantification by considering a generic element, we would like to say which properties are true of the specimen of A. An answer is that the specimen enjoys all the properties P that are true of most As, "most" being defined by the appropriate measure on the type A — and a percentage depending on A and on P. Given that it is never the case that both a property and its negation are true of most of the A, there should not be any contradiction. Because of the duality between properties and individuals, it would be pleasant to also have a measure on the set of predicates that are functions to truth-values. This is mathematically possible, since whenever a set is endowed with a measure, the set of functions from this set to another set can be endowed with a measure as well, as some construction by Kolmogorov shows.

Regarding scalar functions like *height*, *weight*, etc. we would prefer to replace the function by a relation, and to allow the generic element to have a full interval of values rather than a single value. Hence the height of the specimen can be any value in a given interval — an example of such a situation is provided by baby weight and height charts. The idea is that a function "*tall*" means taller than common values, the common values being the ones associated to the specimen, the ones of "most of" the individuals in the relevant class. But it might be trickier than that. In some species, males are taller than females. There should be no problem to have a height chart for males, one for females, but what about the size of a specimen of this species without specifying its sex? One could say that the interval goes from the minimum of the female height interval to the maximum of the male height, but possibly the first value has to be increased and the second one decreased. This is just an intuition that requires further study, with some inspiration from the work of Egré and Klinedinst (2011) and Bale (2011). Our proposal is possibly quite close to the PhD thesis of Solt (2009) that we discovered very recently. As one can see, this is just a proposal, and we are far from a neat answer.

**7.2. Interactive models, proofs and refutations**

Given the complexity of the definition of truth for a sentence involving a vague quantifier, one may adopt a more pragmatic answer: what are the situations, contexts in which such a "most of" statement can be asserted? On the formal side, this leads us to think about the other side of logic, namely proofs. Of course we do not have a complete sets of rules, but still we know some correct ways of reasoning with such a notion. Firstly, when a property is true of *all* individuals, it is true of *"most of"* them, although because of Gricean maxims we do not say *"most of"* when we can say *"all"*. We also know that when a property P is true of *"most"* As there is an A satisfying P — with some precaution about empty models and conditionals with an ironic reading: *"if Brits do not like France, I am the pope"*.

Finding rules for *"most"* can be split into two already difficult questions:
1. Defining rules for proportional quantifiers ("*the majority of*", "*more than 30%*", etc.).
2. Adapting such rules for an undetermined large or small proportion.

"The majority of" is in some sense self-dual: one asserts it because P holds of more than 50%, but a way to refute it is to find another property Q, which also holds of more than 50% but which does not intersect P. The same holds for X% and (100-X)% quantifiers. This, together with our view of *"most"* as acting on a single type, as universal and existential quantifications, suggests that generalized quantifiers apply to predicate and not to sets of individual. In our view, the fact that most A are P is rather a property of P than a property of the elements in A.

We must admit that, for the time being, we can only propose a direction for further investigations. Our proposal is as frustrating as the current literature with tableau rules: monotonicity principles define rules, "*tertium non datur Xv~X for all X*" is the axiom, but tableau rules are also far from complete: nothing distinguishes two quantifiers with the same monotonicity properties, except the models, if they are allowed to intervene inside the rules, which is not so satisfactory, see e.g. Peters and Westerstahl (2008).

Thus the interpretation of the generic elements is an open question, but at least we have a neat syntax for them and know which formulae we want to interpret and how they are obtained from syntactic parse trees, and there exist fine-grained presentations of the linguistic aspects of quantification like Mari (2011) and Szabolcsi (2010).

## 8. On the debate between semantic minimalism and contextualism

We actually started our reflection on generics from classical examples in the minimalism-contextualism debate. These examples are statements that can be both true and false depending on the class in which the object is considered, which is provided by the context. For instance, if Carlotta is a two-year old girl, depending on her class – her type in our type theoretic framework — the following statement can be both true and false:

> (7) Carlotta is tall.
> (8) My daughter is tall and thin for a 2 year old, but she is following her curve." (http://mom4mom.com)
> (9) My two-year-old can't get his own cup out of the cabinet because he can't reach, …( http://633woman.com)

We noticed that the specimen notion together with the flexibility of second-order typing succeeds in capturing this phenomenon. As said above, entries in the lexicon contain optional λ -terms that encode the ontological relations and in the case of a two-year old girl like Carlotta, she can be viewed as a child, and also as a female human being, as a human being etc. Here are the constants and the useful lexicon entries:

- *float=type for real numbers*
- *height : Πα . (α → float → t) height is a binary predicate*
- *<: float→float→t*
- **Carlotta**
    - Carlotta : 2yoGirl (constant)
    - h : 2yoGirl → human (optional λ -term)
- **tall**
    - Λαλx:α
      ∀{float}λh:float
      ∀{float}λhs:float height{α}(∠{α},hs)∧height{α}(x,h)
      ⇒ hs ≤ h
    - **type of tall**: Πα.α →t

The constant *height* is a relation between members of a type and numbers (float), and numbers are compared with <. The entry for *tall* applies to any type T (second order is quite important here as well) and to a term u of type T. It says that the object u is taller than any possible height of the specimen of this class T.

If we do not use any optional λ-term, we apply the lambda term associated to *tall* to the type *2yoGirl*, and to the constant *Carlotta2yoGirl* we get the reading where *Carlotta* is taller than the maximal height of the *2yoGirl* specimen (think again of baby height charts). This is likely to be interpreted as true.

But if we apply "tall" to the *human* type, we cannot apply the result to the constant *Carlotta:2yoGirl*. But we can firstly apply the *h : 2yoGirl → human* (optional λ -term) to the constant *Carlotta:2yoGirl* and proceed: using the type human since *h(Carlotta)* is of type *human*. We thus obtain the formula meaning that Carlotta is tall as a human being, which is unlikely to be interpreted as true.

The semantic machinery produces every possible reading and the context intervenes as a preference for some optional transformation(s). It should be discussed whether there are one or several natural types for an object. Our model can handle any solution: a single natural type, several privileged types... — quite often, such ontoliged or metaphysical questions spontaneously pop up when dealing with the organization of the concepts in the lexicon.

This case shows a general idea underlying our model: terms represent the computational process for obtaining semantic representations, while types that are flexible are pragmatically inferred from the context. This could be called a type-theoretic viewpoint in the debate between contextualism and semantic minimalism.

## 9. Conclusion

We presented a type-theoretic framework for the logical syntax of lexical and compositional semantics, and we focused on its use for "most" quantifiers, depicted via the corresponding generic elements that we called *specimens*. This strongly relies on the second order lambda calculus, with flexible types, as the right framework for meaning assembly.

Along the way we were able to have a viewpoint on the border between semantics and pragmatics. Formally, *semantics* is carried out by the *terms*, while *pragmatics* provides the *types*: types do not drive the computation but they filter impossible readings and sometimes trigger alternative readings by changing the type, hence the comparison class.

This work has been implemented by Richard Moot as part of a large lexicalised categorical French grammar producing semantic representations as Discourse Representation Structure with lambda-DRT (Moot *et al.* (2011)). The semantic representations are expressed in a multi-sorted first-order logic — instead of multi-sorted higher order as we did here.

Delimiting the syntactic part of semantics - finding the proper logical syntax for compositional semantics - yields more questions than it solves... At least it makes the questions clear: how to interpret the semantic representations, the logical formulae associated with sentences involving "most" quantifiers. Two directions are possible. The standard one, for which we should determine the truth of a statement involving the *specimen*, the generic of "most": this seems to be hardly tractable. The other direction, in the proof theoretical tradition, would be to find rules for asserting and refuting sentences involving *specimens*. Although we shall possibly never find the complete set of rules, we can hope to find convincing subsets by viewing quantification as acting on predicates and not on individuals.

**Thanks** This work owes a lot to Sarah-Jane Conrad (Sprachphilosophie, Universität Bern). Indeed, her talk and our discussions initiated at the *Cerisy*

*Context Conference* on the debate between contextualism and semantic minimalism, lead me to a new connection between logical semantics and type theory, here applied to generic elements. Being rather new to formal semantics, I also would like to thank my "advisors" Claire Beyssade, Francis Corblin, David Nicolas, Philippe Schlenker, and especially Alda Mari who, in addition, organised the Genius workshop that I enjoyed. I also thank the people I work with on related issues, Richard Moot, Vito Michele Abrusci, Nicholas Asher and Zhaohui Luo. I finally thank anonymous reviewers as well as Heather Burnett and Hazel Pearson for their insightful comments.